\documentclass[12pt]{article}
\usepackage{amsmath,amssymb,amsthm,xspace,amscd}

\newtheorem{theo}{Theorem}
\newtheorem{lem}{Lemma}

\newcommand{\bZ}{{\mathbb Z}}
\newcommand{\cA}{{\mathcal A}}
\providecommand{\espan}[1]{\text{span}\left\{ #1\right\}}
\DeclareMathOperator{\ann}{ann}

\title{A class of Locally Nilpotent Commutative Algebras}

\author{Antonio Behn\footnote{Supported by Fondecyt 1070243}, \hspace{.1cm} Alberto Elduque\footnote{Supported by the Spanish Ministerio de
 Educaci\'on y Ciencia
 and FEDER (MTM 2007-67884-C04-02) and by the
Diputaci\'on General de Arag\'on (Grupo de Investigaci\'on de
\'Algebra)}, \hspace{.1cm}  Alicia
Labra\footnote{Part of this research was done while this author was
visiting the University of Zaragoza on a grant from Fondecyt 1070243.}\\
\footnotesize  Departamento de Matem\'aticas,
Fac. de Ciencias, Universidad de Chile\\[-4pt]
\footnotesize  Casilla 653, Santiago - Chile, \\[-4pt]
\footnotesize E-mail: abehn@uchile.cl\\
\footnotesize Departamento de Matem\'aticas e
 Instituto Universitario de Matem\'aticas y Aplicaciones,\\[-4pt]
 \footnotesize  Universidad de Zaragoza, 50009 Zaragoza, Spain,\\[-4pt]
\footnotesize  E-mail: elduque@unizar.es\\
\footnotesize Departamento de Matem\'aticas,  Fac. de Ciencias, Universidad de Chile\\[-4pt]
\footnotesize Casilla 653, Santiago - Chile.\\[-4pt]
\footnotesize E-mail: alimat@uchile.cl\\}

\date{}

\begin{document}

\maketitle

\begin{abstract}
This paper deals with  the variety of commutative nonassociative
algebras satisfying the identity $L_x^3+ \gamma L_{x^3} = 0$, $\gamma \in K$. In \cite{CHL} it is proved that if $\gamma = 0,1$ then any finitely generated algebra
 is nilpotent. Here we generalize this result by proving that if $\gamma \neq  -1$,  then  any such algebra is locally nilpotent. Our results require  characteristic $ \; \neq 2,3$.
\end{abstract}


\section{Introduction}

\bigskip
   Let $\cA$ be the free commutative (nonassociative) algebra on a finite set $\{t_1,\ldots,$

\noindent $t_k\}$ of generators. Then $\cA$ is spanned by the commutative monomials in the $t_i$'s and, therefore, it is naturally graded $\cA=\oplus_{n\geq 1} \cA_n$, with $\cA_1=\espan{t_1,\ldots,t_k}$. In other words $\deg(t_i)=1$ for $i=1,\ldots,k$. (Actually, $\cA$ is $\bZ^k$-graded too, with $t_i$ being homogeneous of degree the $k$-tuple with $1$ in the $i$-th position and $0$'s elsewhere).

\bigskip Let  $  \; L_{x_1}   L_{x_2} \cdots L_{x_n} \; $ be a string of left multiplications by monomials $ \; x_i.$
The {\it length of the string} is $ \; n.$  {\it The total degree
of the string} is $\sum_{i=1}^{n} \deg(x_i).$ The {\it max degree}
is the maximum of $\{\deg(x_1), \deg(x_2), ...,$ $ \deg(x_n) \}.$

\medskip
This paper studies the variety of commutative (nonassociative) algebras
$A$ satisfying  the identity
$$   L_x  L_x  L_x   + \gamma L_{x^3} =0, \quad \gamma \in K $$
whose linearizations are
\begin{eqnarray}
&y(x(xa))+x(y(xa))+x(x(ya))+2\gamma((xy)x)a+\gamma((xx)y)a=0\nonumber\\
\nonumber \\ \text{and}
 \label{f}&f(x,y,z,a) = z(y(xa)) + y(z(xa)) + z(x(ya)) + x(z(ya))+ y(x(za))\nonumber \\ &+ x(y(za))+ 2 \gamma [((xy)z)a+((yz)x)a+((zx)y)a] = 0
\end{eqnarray}

In \cite{CHL} it is proved that if $\gamma = 0,1$ then any finitely
generated algebra is nilpotent. Here we generalize this result by proving that if $\gamma \neq  -1,$  then   any such algebra is locally nilpotent.

 \bigskip
In terms of the operators $L_x$'s and using that $A$ is commutative
we obtain the following identities
\begin{eqnarray}
\label{eq1}&L_xL_xL_y+L_xL_yL_x+L_yL_xL_x+\gamma[2L_{(xy)x}+L_{(xx)y}]=0\\
\nonumber \\
 \label{eq2}&L_x L_y L_z + L_x L_z L_y + L_y L_x L_z + L_y L_z L_x \nonumber \\ &+  L_z L_xL_y   + L_z L_y L_x +2\gamma [L_{(xy)z} +
L_{(yz)x} + L_{(zx)y}] = 0
\end{eqnarray}
\begin{eqnarray}
\label{eq3}&L_{x(xa)}+L_xL_{xa}+L_xL_xL_a+\gamma[ 2L_aL_xL_x+L_aL_{xx}]=0 \\
\nonumber \\
\label{eq4}&L_{y(ax)} + L_y L_{ax} +  L_{x(ay)} + L_x L_{ay} +   L_yL_xL_a\nonumber \\ &   + L_x L_y L_a +2\gamma [L_a L_{xy} +
L_a L_xL_y + L_a L_y L_x] = 0
\end{eqnarray}

\medskip
A string or a linear combinations of strings of the same length is called {\it reducible} if  it is expressible as a linear combination of strings each of which has  the same total degree but of shorter lengths. If $X$ and $Y$ are strings or linear combinations of strings of the same length we will say that $X \equiv Y$ if and only if $X - Y$ is reducible. That is, $X - Y$  is expressible as a linear combination of strings each of which has the same total degree but
length less than the common length of $X$ and $Y.$ When $X \equiv Y$ we say $X$ is equivalent to $Y$.

\medskip
Using the above definition, identities (\ref{eq1}), (\ref{eq2}), (\ref{eq3}) and (\ref{eq4}) become
\begin{eqnarray}
\label{L1}& L_yL_xL_x+ L_xL_xL_y+L_xL_yL_x\equiv 0 \\
 \label{L2}&L_x L_y L_z + L_x L_z L_y + L_y L_x L_z + L_y L_z L_x +  L_z L_x L_y   + L_z L_y L_x
\equiv 0\\
\label{L3}&L_xL_xL_a \equiv -2\gamma L_aL_xL_x\\
\label{L4}& L_yL_xL_a   \equiv  -(L_x L_y L_a  +2\gamma [ L_a
L_xL_y + L_a L_y L_x])
\end{eqnarray}

\section{Reducing string}

 \begin{lem} $L_x  L_{x_1} \cdots L_{x_n}   L_x $ is equivalent to a
 linear combination of strings of the same total length where the $L_x$'s are adjacent.

\begin{proof}
   We proceed by induction on $n$: For $n=1$ we get the string $L_x  L_y  L_x,$.
 Using (\ref{L1}) we have
         $L_x  L_y  L_x \equiv -L_x  L_x  L_y  - L_y  L_x  L_x.$

Assume that the result is true if the distance that the $L_x$'s
are apart is less than $k$. Using (\ref{L4}).
\begin{eqnarray*}
 L_x  L_{x_1} L_{x_2} L_{x_3} \cdots  L_{x_k}   L_x\equiv -[L_{x_1}L_xL_{x_2}+2\gamma(L_{x_2}L_{x_1}L_x+L_{x_2}L_xL_{x_1})]L_{x_3} \cdots  L_{x_k}   L_x
\end{eqnarray*}

 \medskip
Altogether, there are
three strings represented in the above expression. In each of
these three strings, the $L_{x}$'s are less than $k$ apart.

By induction each is equivalent to a linear combination of strings
where the $L_x$'s are adjacent.
 \end{proof}
\end{lem}

\medskip

\begin{lem}   The string  $L_x  L_x  L_y  L_y$   is reducible  if $\gamma \neq \pm \; \frac{ 1}{2}.$

\begin{proof}
\begin{eqnarray*}
 (L_xL_xL_y)L_y&\equiv&-2\gamma(L_yL_xL_x)L_y \text{ using (\ref{L3})}\\
&\equiv&2\gamma(L_xL_yL_x+L_xL_xL_y)L_y\text{ using (\ref{L1})}\\
&\equiv&2\gamma L_x(L_yL_xL_y+L_xL_yL_y)\\
&\equiv&-2\gamma L_x(L_yL_yL_x)\text{ using (\ref{L1})}\\
&\equiv&4\gamma^2 L_x(L_xL_yL_y)\text{ using (\ref{L3}) with $x$ and $y$ interchanged.}
\end{eqnarray*}

Therefore $(4\gamma^2-1)L_xL_xL_yL_y$ is reducible, and since $\gamma\neq \pm\frac12$, $L_x  L_x  L_y  L_y$ is reducible.

\medskip

\end{proof}
\end{lem}

\medskip

\begin{lem}\label{Lx3reducible}   $L_x  L_x  L_x $ is reducible.

\begin{proof}   $L_x  L_x  L_x  = - \gamma L_{x^3} \equiv 0 $ \end{proof}
\end{lem}

\begin{lem}   $L_x L_x  L_{x_1}  L_{x_2} \cdots L_{x_k}   L_y  L_y $  is
reducible if  $\gamma \neq \pm \; \frac{ 1}{2}.$

\begin{proof}  We use induction on $k.$
        If $k = 0, \;  L_x  L_x  L_y  L_y $ is reducible using Lemma 2.
      Now in general, using (\ref{L1}),
\begin{eqnarray*}
 L_x  L_x  L_{x_1}L_{x_2}  \cdots L_{x_k} L_y L_y&\equiv&-2\gamma L_{x_1}L_xL_xL_{x_2}\cdots L_{x_k} L_y L_y\\
&\equiv&0\text{ by induction.}
\end{eqnarray*}
\end{proof}
\end{lem}

\medskip
Let $A$ be the free commutative (but not associative) algebra with
$k$ generators.  Let $\dim[n,k]$ be the dimension of the subspace
of $A$ which is spanned by terms of degree less than $n.$  Thus
$\dim[n,k]$ is the number of distinct
 monomials of $A$ with degree less than $n.$

The proof of the following results depends on the above Lemmas and
it is similar to the one given in \cite{CHL}.

\begin{lem}
Let $A$ be the free commutative (but not associative) algebra with
$k$ generators satisfying the identity $L_x^3 + \gamma L_{x^3} =
0$ with $\gamma\neq\pm\frac12$.   Then any string of total degree
$\geq $ $n\dim[n,k]$ is reducible to  strings whose max degree is
$\geq \; n$ or that have an adjacent pair of identical $L_i's$
\end{lem}

\begin{theo}\label{th:reducing}
Let $A$ be the free commutative (but not associative) algebra with
$k$ generators satisfying identity $L_x^3+ \gamma L_{x^3} = 0$
with $\gamma \neq \pm \; \frac{ 1}{2}.$   Then any string of total
degree $ \geq 2n\dim[n,k]+(n-2)$ is reducible to a linear
combination of strings of maximal degree greater than or equal to
$n.$

\end{theo}

\medskip

\section{Nilpotency}

\medskip

In this section $A$ will be a commutative algebra satisfying
the identity $L_x^3+ \gamma L_{x^3} = 0.$

\smallskip

\noindent We define the function $J(x,y,z)$ by
$J(x,y,z)=(xy)z+(yz)x+(zx)y.$

\begin{lem}\label{Wideal}
 Let $A$ be a commutative algebra over a field of
characteristic $\neq 2 $ that satisfies the identity $L_x^3+
\gamma L_{x^3} = 0 $ with $\gamma \neq 0,-1$. Let $W$ be the
linear subspace of $A$ generated by the elements of the form $x^3$
with $x\in A.$ Then
 $W$ is an ideal of $A.$

 \begin{proof} Recall the first linearization of  $L_x^3+ \gamma L_{x^3} = 0$:
$$y(x(xa))+x(y(xa))+x(x(ya))+2\gamma((xy)x)a+\gamma((xx)y)a=0$$
Replacing $a $ by $x$ and $y$ by $a$ we obtain $ax^3 + x(ax^2) +
x(x(xa)) + 2\gamma x(x(xa)) + \gamma x(x^2a) = 0.$ That is,
\begin{equation}\label{eq:star}
(1 +2\gamma) x(x(xa)) +(1 + \gamma) x(x^2a) + x^3 a = 0.
\end{equation}
On the other hand, we have
\[
J(xa,x,x) =2x(x(xa)) + x^2(xa), \;
\; J(x^2,x,a) = x^3a + x^2(xa) + x(x^2a).
\]
Subtracting both identities we have
\[
x(x^2a) = J(x^2,x,a) -J(xa,x,x) +2x (x(xa)) - x^3a.
\]
Replacing this value in \eqref{eq:star} and reordering we have
$$ (3 + 4  \gamma) x(x(xa)) + (1 +  \gamma) [J(x^2,x,a) -
J(xa,x,x)]-  \gamma x^3 a = 0 \;\; \; (**)$$ Using that $x(x(xa)) +
\gamma x^3a = 0,$ we obtain that
$$-4\gamma (1+\gamma) x^3 a +(1 + \gamma)  [J(x^2,x,a) -
J(xa,x,x)]  = 0$$ Since  $\gamma \neq 0,-1,$ we obtain that
$$ x^3a = \frac{1}{4 \gamma}[J(x^2,x,a) -
J(xa,x,x)] \subseteq W$$ since $J(x,a,a) = \frac{1}{2}[(x+a)^3 -
(x-a)^3] - a^3\subseteq W$ and $ J(x,y,z) = \frac{1}{2} ( J(x+z,
x+ z, y) - J(x-z, x-z, y))\subseteq W.$ Therefore $W$ is an ideal
of $A.$

\end{proof}

\end{lem}

\medskip
If $A$ satisfies the identity $L_x^3
+ \gamma L_{x^3} = 0$ and $\gamma \neq -1$, then  $A$ satisfies
the identity
\begin{equation}\label{x4}
 x(x(xx))  = 0
\end{equation}

Linearizing  completely (\ref{x4}) we get:
\begin{eqnarray}
 0 =g(x,y,z,a) = &(a(x(yz)) +a(y(xz))+a(z(xy)) +x(a(yz)) \nonumber\\
               &+x(y(az))+x(z(ay))+ y(a(xz)) +y(x(az))\nonumber\\
 \label{g}              &+y(z (ax)) +z(a(xy))+z(x(ya))+z(y(ax)).
\end{eqnarray}

 \medskip

\begin{theo}\label{thW2}
Let $A$ be a commutative algebra over a field of characteristic
$\neq 2,3 $ that satisfies the identity $L_x^3+ \gamma L_{x^3} = 0
$ with $\gamma(\gamma^2-1)(4\gamma^2-1) \neq 0$.  Then
 the ideal $W$ of $A$ in Lemma \ref{Wideal} satisfies
$W^2 =0.$
\end{theo}

\begin{proof}
Consider the free commutative algebra $\cA$ on two variables $x$ and $y$. This is $\bZ^2$-graded: $\cA=\oplus_{n,m\geq 1}\cA_{n,m}$. Also, there is a natural order $2$ automorphism $\phi$ which permutes these two variables. The automorphism $\phi$ satisfies $\phi(\cA_{n,m})=\cA_{m,n}$ for any $n,m$. In particular it restricts to a linear order $2$ automorphism of $\cA_{3,3}$. The natural basis for the subspace of fixed elements by $\phi$ (that is, the subspace of elements which are symmetric in the two variables) in $\cA_{3,3}$ is the following set $\mathcal{B}$ consisting of $27$ elements:

\medskip
\noindent
$
 \mathcal{B}=\{(x(x(x(y(yy)))))+(y(y(y(x(xx))))),
 (x(x(y(x(yy)))))+(y(y(x(y(xx))))),\\
 (x(x(y(y(xy)))))+(y(y(x(x(xy))))),
 (x(x((xy)(yy))))+(y(y((xx)(xy)))),\\
 (x(y(x(x(yy)))))+(y(x(y(y(xx))))),
 (x(y(x(y(xy)))))+(y(x(y(x(xy))))),\\
 (x(y(y(x(xy)))))+(y(x(x(y(xy))))),
 (x(y(y(y(xx)))))+(y(x(x(x(yy))))),\\
 (x(y((xx)(yy))))+(y(x((xx)(yy)))),
 (x(y((xy)(xy))))+(y(x((xy)(xy)))),\\
 (x((xx)(y(yy))))+(y((yy)(x(xx)))),
 (x((xy)(x(yy))))+(y((xy)(y(xx)))),\\
 (x((xy)(y(xy))))+(y((xy)(x(xy)))),
 (x((yy)(x(xy))))+(y((xx)(y(xy)))),\\
 (x((yy)(y(xx))))+(y((xx)(x(yy)))),
 ((xx)(x(y(yy))))+((yy)(y(x(xx)))),\\
 ((xx)(y(x(yy))))+((yy)(x(y(xx)))),
 ((xx)(y(y(xy))))+((yy)(x(x(xy)))),\\
 ((xx)((xy)(yy)))+((yy)((xx)(xy))),
 ((xy)(x(x(yy))))+((xy)(y(y(xx)))),\\
 ((xy)(x(y(xy))))+((xy)(y(x(xy)))),
 ((xy)((xx)(yy))),
 ((xy)((xy)(xy))),\\
 ((x(xx))(y(yy))),
 ((x(xy))(x(yy)))+((y(xx))(y(xy))),
 ((x(xy))(y(xy))),\\
 ((x(yy))(y(xx)))\}
$

\smallskip

On the other hand, consider the following family $\mathcal{L}$ of $27$ elements in the subspace of symmetric elements in $\cA_{3,3}$ obtained from specializations of the full linearization $f(x,y,z,a)$ of the identity $L_x^3+ \gamma L_{x^3} = 0 $ as in (\ref{f}):

\medskip
\noindent
$
\mathcal{L}= \{f(y,y,y,(xx)x)+f(x,x,x,(yy)y),    f(x,y,y,(xy)x)+f(y,x,x,(yx)y),\\ f(yy,x,y,xx)+f(xx,y,x,yy),  f(xx,y,y,xy)+f(yy,x,x,yx),  f(xy,x,y,xy),\\f((xx)y,y,y,x)+f((yy)x,x,x,y),  f((xy)x,y,y,x)+f((yx)y,x,x,y),\\  f((yy)x,x,y,x)+f((xx)y,y,x,y),   f((yy)y,x,x,x)+f((xx)x,y,y,y),\\f(xx,yy,y,x)+f(yy,xx,x,y),  f(xy,xy,y,x)+f(yx,yx,x,y),\\  f(xy,yy,x,x)+f(yx,xx,y,y), f(y,y,y,xx)x+f(x,x,x,yy)y, \\ f(x,y,y,xy)x+f(y,x,x,yx)y,  f(x,x,y,yy)x+f(y,y,x,xx)y,\\ f(xy,y,y,x)x+f(yx,x,x,y)y,  f(yy,x,y,x)x+f(xx,y,x,y)y, \\ f(xx,y,y,y)x+f(yy,x,x,x)y,  f(xy,x,y,y)x+f(yx,y,x,x)y, \\  f(yy,x,x,y)x+f(xx,y,y,x)y, (f(x,y,y,y)x)x+(f(y,x,x,x)y)y, \\ (f(y,y,y,x)x)x+(f(x,x,x,y)y)y,  (f(x,x,y,y)x)y+(f(y,y,x,x)y)x,\\  (f(y,y,x,x)x)y+(f(x,x,y,y)y)x, (xx)f(x,y,y,y)+(yy)f(y,x,x,x), \\ (xx)f(y,y,y,x)+(yy)f(x,x,x,y),  (xy)f(x,x,y,y)+(yx)f(y,y,x,x)\}
$.

\smallskip

Note that the specializations of all the elements in $\mathcal{L}$ by means of elements in a commutative algebra $A$ over a field of characteristic $\ne 2,3$ satisfying the identity $L_x^3+ \gamma L_{x^3} = 0 $ is $0$.

\medskip
Let $M$ be the matrix with rows representing the elements of $\mathcal{L}$ in the basis $\mathcal{B}$.

\medskip
Using SAGE (\cite{SAGE}) and MAGMA (\cite{MAGMA}) we compute the determinant of this matrix:
\[\det(M)=-2^{33}3^4\gamma^5(\gamma - 1)(\gamma+1)^{10}(2\gamma - 1)^3(2\gamma + 1)\]
For those values of $\gamma$ for which the determinant is not
zero, we can invert the matrix and therefore we can write any
element of $\mathcal{B}$ as a linear combination of relations in
$\mathcal{L}$. In particular, when $\det(M)\neq 0$, $x^3y^3=0$ is
an identity in our algebra $A$, so we conclude that $W^2=0$.
\end{proof}

\medskip

\begin{theo}\label{mainth}
Any commutative algebra over a field of
characteristic $\neq 2, 3, $
 satisfying the identity $L_x^3+ \gamma L_{x^3} = 0$ with $\gamma\neq 0,\pm1,\pm\frac12$ is locally nilpotent.

More precisely, any
 commutative algebra generated by $k$ elements over a field of
characteristic $\neq 2, 3, $
 satisfying the identity $L_x^3+ \gamma L_{x^3} = 0$ with $\gamma\neq 0,\pm1,\pm\frac12$  is
nilpotent of index at most $2^{4n \dim[n,k] +
2(n-2)}$, where $n$ is the index of nilpotency of the free commutative algebra on $k$ generators satisfying the identity $x^3=0$.

\begin{proof} Any product of total degree $\geq 2^{4n \dim[n,k]+2(n-2) }$ is expressible as a string of length greater than
$4 n \dim[n,k] + 2(n-2).$

By Theorem 1, any string of total degree  $ \geq 2 n
\dim[n,k]+ (n-2)$ in the  finitely generated commutative
algebra is reducible to a linear combination of strings in
which one of the factors is of degree greater than or equal to $n.$
Passing to the homomorphic image satisfying  the identity  $L_x^3+
\gamma L_{x^3} = 0$ with $\gamma\neq 0,\pm1,\pm\frac12$, this
factor of degree greater than $n$ must lie in $W.$

    If we let the length of the string be twice as long, then
there will be two factors from $W.$  On multiplying these strings
out, the result will be zero because $W^2 = 0.$
   This finishes the proof of Theorem 3.
\end{proof}
\end{theo}

\section{Exceptional Cases}

We now look at the five  cases which arose as exceptions in Theorem \ref{mainth}.

\medskip
{\bf Case $\gamma = 0$ or $1$.} In \cite{CHL} it was proved that in these
cases every finitely generated commutative algebra $A$  satisfying
the identity  $L_x^3 = 0,$ or  $L_x^3+  L_{x^3} = 0$ is nilpotent.

\medskip

{\bf Case $\gamma = -1.$} The identity becomes $L_x^3 -  L_{x^3} =
0$ or $ x(x(xy)) = x^3 y.$ We observe that any associative
algebra satisfies this identity. In particular the algebra of
polynomials in a finite set of variables satisfies  $L_x^3 -  L_{x^3} =
0$, it is finitely generated but not nilpotent. Therefore, Theorem \ref{mainth} cannot be extended to this case.

\bigskip

{\bf Case $\gamma = \frac{1}{2}.$} The identity becomes $L_x^3 +
\frac{1}{2} L_{x^3} = 0$.

\smallskip

In this case, replacing $a$ by $x$ in identity (4) with $\gamma = \frac{1}{2}$
we obtain
$$ L_{x^3} + L_x L_{x^2} + L_{x^3} + \frac{1}{2}[2L_{x^3} + L_x
L_{x^2}] = 0$$ That is, $$ L_{x^3}  + 2L_{x^3} + \frac{3}{2}  L_x
L_{x^2} = 0$$ Therefore, $\frac{3}{2}  L_x L_{x^2} = 0$ and
characteristic not $3$ implies that
\begin{equation}\label{eq:LxLx2}
L_x L_{x^2} = 0.
\end{equation}

Now, replacing $a$ by $x^2$ in identity (4) with $\gamma = \frac{1}{2}$
we obtain the identity $$ L_{x^4} + L_x  L_{x^3} + L_x L_x L_{x^2}
+  L_{x^2} L_x L_x + \frac{1}{2}  L_{x^2}  L_{x^2} = 0 $$
where $x^4 = x^3 x.$

Using
that $ x^4 = 0,  L_{x^3} = - 2  L_{x^3}$ and $ L_x L_{x^2} = 0$ we
have that
\begin{equation}\label{eq:Lx4..}
 -2L_x^4  + L_{x^2} L_x^2 + L_{x^2}L_{x^2} = 0.
 \end{equation}

Multiplying \eqref{eq:Lx4..} to the left by $L_x$ and using \eqref{eq:LxLx2} we obtain
\begin{equation}\label{Lx5}
L_x^5 = 0,
\end{equation}
so $L_x$ is nilpotent, and by linearization of $L_x^5$, using \eqref{eq:LxLx2}, we get
\[
L_{x^2}L_x^4=0.
\]
Linearize \eqref{eq:LxLx2} to obtain $2L_xL_{xy}+L_yL_{x^2}=0$ which, with $y=x^2$ gives
\[
(L_{x^2})^2=4L_x^4,
\]
which transforms, using \eqref{eq:Lx4..}, into
\begin{equation}\label{eq:Lx2Lx2}
L_{x^2}L_x^2=0.
\end{equation}
Finally, with $y=x^2$ in \eqref{L1}, equations \eqref{eq:LxLx2} and \eqref{eq:Lx2Lx2} give
\begin{equation}\label{eq:Lx2x2}
L_{x^2x^2}=0.
\end{equation}
That is, for every $ \;
x \in A$, $x^2 x^2 \in \ann(A)=\{z\in A: zA=0\}$.

\begin{theo}
Any  commutative algebra over a field of
characteristic $\neq 2, 3, $ satisfying the identity $L_x^3+
\frac{1}{2}L_{x^3} = 0$, is locally nilpotent.
\end{theo}

\begin{proof}
Let $A$ be a commutative algebra over a field of
characteristic $\neq 2, 3, $ satisfying the identity $L_x^3+
\frac{1}{2}L_{x^3} = 0$. It is enough to prove that $\tilde A=A/\ann(A)$ is locally nilpotent.

But according to \eqref{eq:Lx2x2}, $\tilde A$ satisfies the equation $x^2x^2=0$, and hence, by linearization, the equation $(xy)x^2=0$. On the other hand, equation \eqref{eq:LxLx2} shows that $x(yx^2)=0$. In particular, $\tilde A$ is a Jordan algebra. Let $W$ be the ideal spanned by its cubes (Lemma \ref{Wideal}). Since $(x^3)^2=L_{x^3}(x^3)= 2x(x(xx^3)=0$, and $L_{x^3}^2=4L_x^6=0=L_{(x^3)^2}$, it turns out that the cubes are absolute zero divisors of the Jordan algebra $\tilde A$. Thus $W$ is a locally nilpotent ideal of $\tilde A$, because of Zel'manov's Local Nilpotence Theorem (see \cite[p.~1004]{McCrimmonbook} and \cite{Zelmanov}). Also the algebra $\tilde A/W$ satisfies the identity $x^3=0$, so it is locally nilpotent (see \cite[p.~114]{ZSSS}). Therefore the Jordan algebra $\tilde A$ is locally nilpotent (\cite[Chapter 4, Lemma 7]{ZSSS}), and so is $A$, as required.
\end{proof}

\bigskip

{\bf Case $\gamma = -\frac{1}{2}.$} The identity becomes $L_x^3 =
\frac{1}{2} L_{x^3}$.

\smallskip

In this case, our algebra $A$ satisfies $x^3x=0$ too. Also, $x^3x^2=2x(x(xx^2))=0$, and with $y=x^2$ in \eqref{eq1} it follows that $L_x^5=0$. Now, identities \eqref{L1}, \eqref{L3} and \eqref{L4} show that
\begin{gather}\label{eq:Lx2Ly}
L_x^2L_y\equiv L_yL_x^2,\\
L_xL_yL_x\equiv -2L_yL_x^2,
\end{gather}
and these, together with Lemma \ref{Lx3reducible}, immediately imply the following result

\begin{lem} Any string of the form $L_x  L_{x_1} \cdots L_{x_n}   L_x $ is equivalent to a linear combination of strings of the same total length which end up in $L_x^2$. Also, any string containing three equal elements $L_x$ is reducible.
\end{lem}

There is the following counterpart to Theorem \ref{th:reducing}

\begin{theo}
Let $A$ be the free commutative (but not associative) algebra with
$k$ generators satisfying identity $L_x^3-\frac{1}{2} L_{x^3} = 0$.
Then any string of total
degree $ \geq 2n\dim[n,k]+(n-2)$ is reducible to a linear
combination of strings of maximal degree greater than or equal to
$n$.
\end{theo}
\begin{proof}
If the total degree of a string $L_{x_1}\cdots L_{x_m}$ is $\geq 2n\dim[n,k]$ and its maximal degree is $<n$, then its length $m$ is greater than $2\dim[n,k]$, and therefore there are three monomials $x_i$ which are equal. By the previous Lemma this is reducible. The process can be continued until the maximal degree be $\geq n$, as required.
\end{proof}

And there is too the counterpart to Theorem \ref{thW2}:

\begin{theo}
Let $A$ be a commutative algebra over a field of characteristic
$\neq 2,3 $ that satisfies the identity $L_x^3=\frac{1}{2} L_{x^3}$.  Then
the ideal $W$ of $A$ in Lemma \ref{Wideal} satisfies
$W^2 =0$.
\end{theo}

\begin{proof}
Consider again the free commutative algebra $\cA$ on two variables $x$ and $y$ and the subspace of symmetric elements in $\cA_{3,3}$ spanned by the linearly independent subset (consisting of $19$ elements):

\medskip
\noindent
$
 \mathcal{B'}=\{x(x(x(y(yy)))) + y(y(y(x(xx)))),
 x(x(y(x(yy)))) + y(y(x(y(xx)))),\\
  x(x(y(y(xy)))) + y(y(x(x(xy)))),
 x(x((xy)(yy))) + y(y((xx)(xy))),\\
  x(y(x(x(yy)))) + y(x(y(y(xx)))),
 x(y(x(y(xy)))) + y(x(y(x(xy)))),\\
  x(y(y(x(xy)))) + y(x(x(y(xy)))),
 x(y(y(y(xx)))) + y(x(x(x(yy)))),\\
  x(y((xx)(yy))) + y(x((xx)(yy))),
 x(y((xy)(xy))) + y(x((xy)(xy))),\\
  x((xx)(y(yy))) + y((yy)(x(xx))),
 x((xy)(x(yy))) + y((xy)(y(xx))),\\
  x((xy)(y(xy))) + y((xy)(x(xy))),
 x((yy)(x(xy))) + y((xx)(y(xy))),\\
  x((yy)(y(xx))) + y((xx)(x(yy))),
  (x(xy))(x(yy)) + (y(xx))(y(xy)),\\
 (x(xy))(y(xy)),
 (x(yy))(y(xx)), (x(xx))(y(yy))\}
$.

\smallskip

On the other hand, consider the following family $\mathcal{L'}$ of $19$ elements in the subspace of symmetric elements in $\cA_{3,3}$ obtained from specializations of the full linearization $f(x,y,z,a)$ of the identity $L_x^3-\frac{1}{2}L_{x^3} = 0 $ as in (\ref{f}):

\medskip\noindent
$ \mathcal{L'}= \{f(y,y,y,(xx)x) + f(x,x,x,(yy)y),  f(x,y,y,(xx)y)
+ f(y,x,x,(yy)x),\\ f(x,y,y,(xy)x) + f(y,x,x,(yx)y),
f((xx)y,y,y,x) + f((yy)x,x,x,y),\\  f((xy)x,y,y,x) +
f((yx)y,x,x,y),  f((yy)x,x,y,x) + f((xx)y,y,x,y),\\
f((yx)y,x,y,x) + f((xy)x,y,x,y),  f((yy)y,x,x,x) +
f((xx)x,y,y,y),\\  f(y,y,y,xx)x + f(x,x,x,yy)y,  f(x,y,y,xy)x +
f(y,x,x,yx)y,\\  f(x,x,y,yy)x + f(y,y,x,xx)y, f(xy,y,y,x)x +
f(yx,x,x,y)y,\\  f(yy,x,y,x)x + f(xx,y,x,y)y,  f(xy,x,y,y)x +
f(yx,y,x,x)y,\\  f(yy,x,x,y)x + f(xx,y,y,x)y,  (f(x,y,y,y)x)x +
(f(y,x,x,x)y)y,\\  (f(y,y,y,x)x)x + (f(x,x,x,y)y)y,
(f(x,x,y,y)x)y + (f(y,y,x,x)y)x,\\  f(xx,y,y,y)x + f(yy,x,x,x)y\}
$.

\medskip
Let $M'$ be the matrix with rows representing the elements of $\mathcal{L'}$ in the basis $\mathcal{B'}$. Using again SAGE (\cite{SAGE})
and MAGMA (\cite{MAGMA}) we compute the determinant of this matrix:
\[
\det(M')=2^{14}3^4.
\]
In particular, this shows that in any commutative algebra satisfying the identity $L_x^3=\frac12 L_{x^3}$, $x^3y^3=0$ for any $x,y$, so that $W^2=0$.
\end{proof}

\medskip

Finally the same arguments as for Theorem \ref{mainth} settle this exceptional case:

\begin{theo}
Any commutative algebra over a field of
characteristic $\neq 2, 3, $
 satisfying the identity $L_x^3=\frac12 L_{x^3} $ is locally nilpotent.

More precisely, any
 commutative algebra generated by $k$ elements over a field of
characteristic $\neq 2, 3, $
 satisfying the identity $L_x^3=\frac12 L_{x^3}$  is
nilpotent of index at most $2^{4n \dim[n,k]}$, where $n$ is the index of nilpotency of the free commutative algebra on $k$ generators satisfying the identity $x^3=0$.
\end{theo}

\end{document}